\documentclass[12pt]{amsart}

\allowdisplaybreaks
\usepackage{amssymb}
\usepackage{graphicx}
\usepackage{fullpage}
\usepackage[pdfencoding=auto,colorlinks,linkcolor=,citecolor=]{hyperref}
\usepackage{cite}
\usepackage{dsfont}
\usepackage{eucal}
\usepackage[all,cmtip]{xy}

\newcommand{\Br}{\operatorname{\sf Br}}
\newcommand{\End}{\operatorname{\sf End}}
\newcommand{\Gal}{\operatorname{\sf Gal}}
\newcommand{\GL}{\operatorname{\sf GL}}
\newcommand{\Hom}{\operatorname{\sf Hom}}
\newcommand{\ind}{\operatorname{\sf ind}}
\newcommand{\op}{{\mathsf{op}}}
\newcommand{\res}{\operatorname{\sf res}}
\newcommand{\bQ}{\mathbb{Q}}

\numberwithin{equation}{section}  

\theoremstyle{plain}
\newtheorem{lemma}[equation]{Lemma}
\newtheorem{theorem}[equation]{Theorem}
\newtheorem{proposition}[equation]{Proposition}
\newtheorem{corollary}[equation]{Corollary}

\theoremstyle{definition}

\newtheorem{definition}[equation]{Definition}

\theoremstyle{remark}

\author{David J. Benson} 
\address{Institute of Mathematics \\ 
University of Aberdeen \\ 
Aberdeen AB24 3UE \\ 
United Kingdom}
\email{d.j.benson@abdn.ac.uk}



\title{Matrices for finite group representations that respect Galois automorphisms}

\begin{document}

\begin{abstract}
 We are given a finite group $H$, an automorphism $\tau$ of $H$ of
 order $r$, a Galois
 extension $L/K$ of fields of characteristic zero with cyclic Galois
 group $\langle\sigma\rangle$ of order $r$, and an
 absolutely irreducible representation $\rho\colon H\to\GL(n,L)$ such
 that the action of $\tau$ on the character of $\rho$ is the same as
 the action of $\sigma$. Then the following are equivalent.

 $\bullet$ $\rho$ is equivalent to a
 representation $\rho'\colon H\to\GL(n,L)$ such that the action of
 $\sigma$ on the entries of the matrices corresponds to the action of
 $\tau$ on $H$, and

 $\bullet$ the induced representation
 $\ind_{H,H\rtimes\langle\tau\rangle}(\rho)$ has Schur index one; that
 is, it is similar to a representation over $K$.

 As examples, we
 discuss a three dimensional irreducible representation of $A_5$ over
 $\bQ[\sqrt5]$ and a four dimensional irreducible representation of
 the double cover of $A_7$ over $\bQ[\sqrt{-7}]$.
\end{abstract}

\maketitle

\section{Introduction}

This paper begins with the following question, suggested to
the author by Richard Parker. The alternating
group $A_5$ has a three dimensional representation over the
field $\bQ[\sqrt5]$ which induces up to the symmetric group
$S_5$ to give a six dimensional irreducible that can be written over $\bQ$.
Given an involution in $S_5$ that is not in $A_5$, is it possible to
write down a $3\times 3$ matrix representation of $A_5$ such that the Galois
automorphism of $\bQ[\sqrt5]$ acts on matrices in the same way as the involution
acts on $A_5$ by conjugation?

More generally, we are given a finite group $H$, an automorphism
$\tau$ of order $r$, a Galois extension $L/K$ of fields of characteristic zero
with cyclic Galois group
$\Gal(L/K)=\langle\sigma\rangle$ of order $[L:K]=r$, and an absolutely irreducible
representation $\rho\colon H\to\GL(n,L)$.
We assume that the action of $\tau$ on the character of the representation
$\rho$ is the same as the action of $\sigma$.
Then the question is whether it is possible to conjugate to a representation
$\rho'\colon H\to\GL(n,L)$ with the property that the Galois automorphism $\sigma$
acts on matrices in the same way as $\tau$ acts on $H$. In other words,
we are asking whether the following diagram can be made to commute.
\[ \xymatrix{
    H\ar[r]^(0.4){\rho'}\ar[d]_\tau&\GL(n,L)\ar[d]^\sigma\\
    H\ar[r]^(0.4){\rho'}&\GL(n,L)} \]

We answer this using the invariant $\lambda(\rho)$ in the relative Brauer group
\[ \Br(L/K)=H^2(\langle\sigma\rangle,L^\times)\cong
K^\times/N_{L/K}(L^\times) \]
that defines the division algebra associated to the 
representation obtained by inducing to the semidirect product.

\begin{theorem}\label{th:main}
  Let $\rho\colon H \to \GL(n,L)$ be as above. Then there is an
  invariant $\lambda(\rho)\in K^\times/N_{L/K}(L^\times)$ such that
  the following are equivalent.
\begin{enumerate}
\item $\lambda(\rho)=1$,
\item There is a conjugate $\rho'$ of $\rho$ making the diagram
above commute,
\item If $G$ is the semidirect product $H\rtimes \langle \tau\rangle$
then the induced representation $\ind_{H,G}(\rho)$ has Schur index
equal to one; in other words, it can be written over $K$.
\end{enumerate}
More generally, the order of $\lambda(\rho)$ in
$K^\times/N_{L/K}(L^\times)$
is equal to the Schur index of the induced representation, and the
associated division algebra is the one determined by $\lambda(\rho)$.
\end{theorem}

The equivalence of (1) and (2) is proved in Section~\ref{se:Y}.
The equivalence of (1) and (3) is more standard,
see for example Turull~\cite{Turull:2000a}, and is proved in
Section~\ref{se:ind}. Combining these gives the more interesting
statement of the equivalence of (2) and (3).
We end with some examples. In the case of the three dimensional
representations of $A_5$, we have $\lambda(\rho)=1$,
and we write down explicit matrices for $\rho'$, though they're not
very pleasant. In the case of the
four dimensional irreducible representations of $2A_7$, we have
$\lambda(\rho)=-2$, which is not a norm  from $\bQ[\sqrt{-7}]$,
and the division ring associated to the induced representation is the
quaternion algebra with symbol $(-2,-7)_\bQ$.\bigskip

\noindent
{\bf Acknowledgement.} I would like to thank Richard Parker for
suggesting this problem and the examples, and Alexandre Turull for
some helpful comments.

\section{The matrix $X$}\label{se:X}

Consider the composite $\sigma\circ\rho\circ\tau^{-1}$:
\[ H \xrightarrow{\tau^{-1}} H \xrightarrow{\rho}\GL(n,L) \xrightarrow{\sigma}\GL(n,L). \]
This representation is equivalent to $\rho$, and so there exists a
matrix $X$, well defined up to scalars in $L^\times$, such that conjugation by $X$ takes $\rho$ to $\sigma\circ\rho\circ\tau^{-1}$.
Write $c_X$ for conjugation by $X$, so that $c_X(A)=XAX^{-1}$. Then we have
\begin{equation}\label{eq:1}
  \sigma\circ\rho\circ\tau^{-1} = c_X \circ\rho.
\end{equation}
By abuse of notation, we shall also write $\sigma$ for the automorphism
of $\GL(n,L)$ given by applying $\sigma$ to each of its entries. Then
$c_{\sigma(X)}(\sigma(A))=\sigma(X)\sigma(A)\sigma(X)^{-1}=\sigma(XAX^{-1})$,
so we have
\[ c_{\sigma(X)}\circ\sigma = \sigma\circ c_X. \]
So equation~\eqref{eq:1} gives
\begin{align*}
\sigma^2\circ\rho\circ\tau^{-2}&=\sigma\circ c_X\circ\rho\circ\tau^{-1}\\
&=c_{\sigma(X)}\circ\sigma\circ\rho\circ\tau^{-1}\\
&=c_{\sigma(X)}\circ c_X \circ\rho\\
&=c_{\sigma(X).X}\circ\rho.
\end{align*}
Continuing this way, for any $i>0$ we have
\[ \sigma^i\circ\rho\circ\tau^{-i} = c_{\sigma^{i-1}(X)\cdots\sigma(X).X}\circ\rho. \]
Taking $i=r$, we have $\sigma^r=1$ and $\tau^r=1$, so
\begin{equation}\label{eq:NX}
  \rho=c_{\sigma^{r-1}(X)\cdots\sigma(X).X}\circ\rho.
\end{equation}

\begin{definition}
  If $A$ is an $n\times n$ matrix over $L$, we define the norm of $A$ to be
  \[ N_{L/K}(A)=\sigma^{r-1}(A)\cdots\sigma(A).A \]
  as an $n\times n$ matrix over $K$.
\end{definition}

Equation~\eqref{eq:NX} now reads
\[ \rho=c_{N_{L/K}(X)}\circ\rho. \]
By Schur's lemma, it follows that the matrix
$N_{L/K}(X)$ is a scalar multiple of the identity,
\[ N_{L/K}(X)=\lambda I. \]
Applying $\sigma$ and rotating the terms on the left, we see that
$\lambda=\sigma(\lambda)$, so that $\lambda\in K^\times$. If we replace
$X$ by a scalar multiple $\mu X$, then the scalar $\lambda$ gets multiplied
by $\sigma^{r-1}(\mu)\cdots \sigma(\mu)\mu$, which is the norm $N_{L/K}(\mu)$.
Thus the scalar $\lambda$ is well defined only up to norms of elements in $L^\times$.
We define it to be the $\lambda$-invariant of $\rho$:
\[ \lambda(\rho)\in K^\times/N_{L/K}(L^\times). \]
Thus $\lambda(\rho)=1$ if and only if $X$ can be replaced by a multiple of $X$ to
make $N_{L/K}(X)=I$.

\section{The matrix $Y$}\label{se:Y}

The goal is to find a matrix $Y$ conjugating $\rho$ to a representation
$\rho'$ such that $\sigma\circ\rho'\circ\tau^{-1}=\rho'$. Thus we wish $Y$
to satisfy
\[ \sigma\circ c_Y\circ \rho\circ\tau^{-1}=c_Y\circ\rho. \]
We rewrite this in stages:
\begin{align*}
  c_{\sigma(Y)}\circ\sigma\circ\rho\circ\tau^{-1}&=c_Y\circ\rho\\
 \sigma\circ\rho\circ\tau^{-1}&=c_{\sigma(Y)^{-1}}\circ c_Y\circ\rho\\
 c_X\circ\rho&=c_{\sigma(Y)^{-1}Y}\circ\rho.
\end{align*}
Again applying Schur's lemma, $\sigma(Y)^{-1}Y$ is then forced to be a multiple of $X$.
Since $N_{L/K}(\sigma(Y)^{-1}Y)=I$,
it follows that if there is such a $Y$ then $\lambda(\rho)$ is the identity element of
$K^\times/N_{L/K}(L^\times)$. This proves one direction of Theorem~\ref{th:main}. The
other direction is now an immediate consequence of the version of Hilbert's Theorem 90
given in Chapter~X, Proposition~3 of Serre~\cite{Serre:1979a}:

\begin{theorem}\label{th:90}
  Let $L/K$ be a finite Galois extension with Galois
  group $\Gal(L/K)$. Then $H^1(\Gal(L/K),\GL(n,L))=0$.\qed
\end{theorem}

\begin{corollary}
Let $L/K$ be a Galois extension with cyclic Galois group
$\Gal(L/K)=\langle\sigma\rangle$ of order $r$.
If a matrix $X\in\GL(n,L)$ satisfies $N_{L/K}(X)=I$ then
  there is a matrix $Y$ such that $\sigma(Y)^{-1}Y=X$.
\end{corollary}
\begin{proof}
  This is the case of a cyclic Galois group of Theorem~\ref{th:90}.
\end{proof}

This completes the proof of the equivalence of (1)
and (2) in Theorem~\ref{th:main}.

\section{The induced representation}\label{se:ind}

Let $G=H\rtimes\langle\tau\rangle$, so that for $h\in H$ we have
$\tau(h)=\tau h\tau^{-1}$ in $G$. Then the induced representation
$\ind_H^G(\rho)$ is an $LG$-module with character values in $K$, but
cannot necessarily be written as an extension to $L$ of a $KG$-module.
So we restrict the coefficients to $K$ and examine the endomorphism
ring.

\begin{lemma}\label{le:r^2}
  $\End_{KG}(\ind_{H,G}(\rho|_K))$ has dimension $r^2$ over $K$.
\end{lemma}
\begin{proof}
  The representation $\rho|_K$ is an irreducible $KH$-module,
whose extension to $L$ decomposes as the sum of the Galois conjugates
of $\rho$, so $\End_{KH}(\rho|_K)$ is $r$ dimensional over $K$.
For the induced representation
$\ind_{H,G}(\rho|_K)=\ind_{H,G}(\rho)|_K$, as vector spaces we then have
\begin{equation*} \End_{KG}(\ind_{H,G}(\rho|_K))\cong
  \Hom_{KH}(\rho|_K,\res_{G,H}\ind_{H,G}(\rho|_K))\cong
  r.\End_{KH}(\rho|_K). \qedhere
\end{equation*}
\end{proof}

\begin{proposition}\label{pr:End}
The algebra $\End_{KG}(\ind_{H,G}(\rho|_K))$ is a crossed product
algebra, central simple over $K$, with generators $m_\lambda$ for
$\lambda\in L$ and an element 
$\xi$, satisfying
\[ m_\lambda + m_{\lambda'}=m_{\lambda+\lambda'},\qquad
m_\lambda m_{\lambda'}=m_{\lambda\lambda'},\qquad
m_\lambda\circ\xi=\xi\circ m_{\sigma(\lambda)},\qquad
\xi^r=m_{\lambda(\rho)}. \]
\end{proposition}
\begin{proof}
We can write the representation $\ind_{H,G}(\rho|_K)$ in terms of
matrices as follows.
\[ g \mapsto \begin{pmatrix}\rho(g)|_K\\&\sigma\rho\tau^{-1}(g)|_K\\
&&\ddots\\ &&&\sigma^{-1}\rho\tau(g)|_K
\end{pmatrix},\qquad
\tau\mapsto \begin{pmatrix}&&&I\\
I\\&\ddots\\&&I
\end{pmatrix}\circ\sigma \]
It is easy to check that the following are endomorphisms of this representation.
\[ m_\lambda=\begin{pmatrix}\lambda I\\
&\sigma(\lambda) I\\&&\ddots\\
&&&\sigma^{-1}(\lambda) I \end{pmatrix},\qquad
\xi =\begin{pmatrix}
&&&&\sigma^{-1}(X)\\
X\\&\sigma(X)\\
&&\ddots\\
&&&\sigma^{-2}(X)
\end{pmatrix}\]
with $\lambda\in L$ and $X$ as in Section~\ref{se:X}. Since these
generate an algebra of dimension $r^2$ over $K$, by Lemma~\ref{le:r^2} they generate
the algebra $\End_{KG}(\ind_{H,G}(\rho|_K))$. The given relations are easy
to check, and present an algebra which is easy to see has dimension at
most $r^2$, and therefore no further relations are necessary.
\end{proof}

\begin{corollary}
The Schur index of the induced representation $\ind_{H,G}(\rho)$ is
equal to the order of $\lambda(\rho)$ as an element of $K^\times/N_{L/K}(L^\times)$.
In particular, the Schur index is one if and only if $\lambda(\rho)=1$
as an element of $K^\times/N_{L/K}(L^\times)$.
\end{corollary}
\begin{proof}
This follows from the structure of the central simple algebra
$\End_{KG}(\rho|_K)$ given in Proposition~\ref{pr:End}, using the
theory of cyclic crossed product algebras, as developed for example in Section~15.1 of
Pierce~\cite{Pierce:1982a}, particularly Proposition~b of that section.
\end{proof}

This completes the proof of the equivalence of (1) and (3) in
Theorem~\ref{th:main}. In particular, it shows that $\lambda(\rho)$
can only involve primes dividing the order of $G$.
  
\section{Examples}

Our first example is a three dimensional representation of $A_5$.
There are two algebraically conjugate three dimensional irreducible
representations of $A_5$ over
$\bQ[\sqrt5]$ swapped by an outer automorphism of $A_5$, and giving
a six dimensional representation of the symmetric group $S_5$ over $\bQ$.

Setting $\alpha=\frac{1+\sqrt5}{2}$, $\bar\alpha=\frac{1-\sqrt5}{2}$,
we can write the action of the generators
on one of these three dimensional representations as follows.
\[ (12)(34)\mapsto \begin{pmatrix}-1\ &0&0\\0&0&1\\0&\ 1\ &0\end{pmatrix}
  \qquad
  (153)\mapsto\begin{pmatrix}-1&1&\alpha\\
    \alpha&0&-\alpha\\
  -\alpha&0&1\end{pmatrix}\]
Taking this for $\rho$, we find a matrix $X$ conjugating this
to $\sigma\circ\rho\circ\tau^{-1}$ where $\sigma$ is the field automorphism
and $\tau$ is conjugation by $(12)$. Using the fact that if $a=(12)(34)$ and $b=(153)$ then
$ab^2abab^2=(253)$,  we find that
\[ X=\begin{pmatrix}1&-\bar\alpha\ &\bar\alpha\\
    -\bar\alpha\ &1&-\bar\alpha\ \\\bar\alpha&-\bar\alpha\ &1\end{pmatrix}\]
We compute that $\sigma(X).X$ is minus the identity. Now $-1$ is
in the image of $N_{\bQ[\sqrt5],\bQ}$, namely we have $(2-\sqrt5)(2+\sqrt5)=-1$.
So we replace $X$ by $(2-\sqrt5)X$ to achieve $\sigma(X).X=I$.
Having done this, by Hilbert 90 there exists $Y$ with $\sigma(Y)^{-1}.Y=X$.
Such a $Y$ conjugates $\rho$ to the desired form. For example we can take
\[ Y = \begin{pmatrix}
    1-2\sqrt5&3-2\sqrt5&-3+2\sqrt5\\
    3-2\sqrt5&1-2\sqrt5&3-2\sqrt5\\
    -3+2\sqrt5&3-2\sqrt5&1-2\sqrt5
  \end{pmatrix}. \]
Thus we end up with the representation
\[ (12)(34) \mapsto\begin{pmatrix}-1\ &0&0\\0&0&1\\0&\ 1\ &0\end{pmatrix}
  \qquad
  (153)\mapsto\frac{1}{40}\begin{pmatrix}
    10-4\sqrt5&-5+19\sqrt5&25-9\sqrt5\\
    -10-4\sqrt5&25+9\sqrt5&-5-19\sqrt5\\
    -50&35-5\sqrt5&-35-5\sqrt5
  \end{pmatrix} \]
\[ (253)\mapsto
\frac{1}{40}\begin{pmatrix}
    10+4\sqrt5&-5-19\sqrt5&25+9\sqrt5\\
    -10+4\sqrt5&25-9\sqrt5&-5+19\sqrt5\\
    -50&35+5\sqrt5&-35+5\sqrt5
  \end{pmatrix}. \]
Denoting these matrices by $a$, $b$ and $c$, it is routine to check that
$a^2=b^3=(ab)^5=1$, $a^2=c^3=(ac)^5=1$,  and $c=\sigma(b)=ab^2abab^2$.

More generally, if $H$ is an alternating group $A_n$ and $G$ is the
corresponding symmetric group $S_n$ then all irreducible
representations of $G$ are rational and so the invariant
$\lambda(\rho)$ is equal to one for any irreducible character of $H$
that is not rational. So an appropriate matrix $Y$ may always be
found in this case.\bigskip

Our second example is one with $\lambda(\rho)\ne 1$. Let $H$ be the
group $2A_7$, namely a
non-trivial central extension of $A_7$ by a cyclic group of order two.
Let $\tau$ be an automorphism of $H$ of order two, lifting the action
of a transposition in $S_7$ on $H$, and let $G$ be the semidirect
product $H\rtimes\langle\tau\rangle$. Then $H$ has two Galois conjugate
irreducible representations of dimension four over
$\bQ[\sqrt{-7}]$. Let $\rho$ be one of them. The induced
representation is eight dimensional over $\bQ[\sqrt{-7}]$. Restricting
coefficients to $\bQ$ produces 
a $16$ dimensional rational representation whose endomorphism algebra $E$
is a quaternion algebra. Thus the induced
representation can be written as a four dimensional representation
over $E^\op\cong E$. This endomorphism algebra was computed by
Turull~\cite{Turull:1992a} in general for the double covers of
symmetric groups. In this case, by Corollary~5.7 of that paper, the algebra
$E$ is generated over $\bQ$ by elements $u$ and $v$ satisfying $u^2=-2$, $v^2=-7$ and
$uv=-vu$.
Thus the invariant $\lambda(\rho)$ is equal to $-2$ as an element of
$\bQ^\times/N_{\bQ[\sqrt{-7}],\bQ}(\bQ[\sqrt{-7}]^\times)$
in this case.

\bibliographystyle{amsplain}
\bibliography{../repcoh}

\newcommand{\noopsort}[1]{}
\providecommand{\bysame}{\leavevmode\hbox to3em{\hrulefill}\thinspace}
\providecommand{\MR}{\relax\ifhmode\unskip\space\fi MR }
\providecommand{\MRhref}[2]{%
  \href{http://www.ams.org/mathscinet-getitem?mr=#1}{#2}
}
\providecommand{\href}[2]{#2}
\begin{thebibliography}{1}

\bibitem{Pierce:1982a}
R.~S. Pierce, \emph{Associative algebras}, Graduate Texts in Mathematics,
  vol.~88, Springer-Verlag, Ber\-lin/New York, 1982.

\bibitem{Serre:1979a}
J.-P. Serre, \emph{{Local fields}}, Graduate Texts in Mathematics, vol.~67,
  Springer-Verlag, Ber\-lin/New York, 1979.

\bibitem{Turull:1992a}
A.~Turull, \emph{{The Schur index of projective characters of symmetric and
  alternating groups}}, Ann.\ of Math. \textbf{135} (1992), 91--124.

\bibitem{Turull:2000a}
\bysame, \emph{{Clifford theory for cyclic quotient groups}}, J.~Algebra
  \textbf{227} (2000), 133--158.

\end{thebibliography}

\end{document}